\documentclass[conference,a4paper]{IEEEtran}
\IEEEoverridecommandlockouts
\usepackage{array}
\usepackage{epsfig}
\usepackage{enumitem}
\usepackage{graphics}
\usepackage{graphicx}
\usepackage[english]{babel}
\usepackage[latin1]{inputenc}
\usepackage{amsmath}
\usepackage{amssymb}
\usepackage{booktabs}
\usepackage{floatrow}
\usepackage{subfigure}
\usepackage{float}
\usepackage{bm}
\usepackage{cite}
\usepackage{cases}
\usepackage{color,soul}
\usepackage{comment}
\usepackage{multicol}
\usepackage{multirow}
\usepackage{balance}
\usepackage{acronym}
\usepackage{chemformula}
\usepackage{eurosym}
\usepackage{tikz}

\usepackage{siunitx} 
\DeclareSIUnit{\atm}{atm}
\DeclareSIUnit{\kWh}{kWh}
\DeclareSIUnit{\Ah}{Ah}
\DeclareSIUnit{\eur}{\mbox{\text{\euro}}}
\DeclareSIUnit{\kW}{kW}
\DeclareSIUnit{\MW}{MW}
\DeclareSIUnit{\kwm}{kW/m^2}
\DeclareSIUnit{\m}{m}
\DeclareSIUnit{\h}{h}

\restylefloat{table}
\newcommand{\modifPUMA}[1]{\textcolor{black}{#1}}

\acrodef{rfc}[RFC]{Reverse Fuel Cell}
\acrodef{soec}[SOEC]{Solid-Oxide Electrolyzer Cell}
\acrodef{sofc}[SOFC]{Solid-Oxide Fuel Cell}
\acrodef{tsoec}[t-SOEC]{transition to Solid-Oxide Electrolyzer Cell}
\acrodef{tsofc}[t-SOFC]{transition to Solid-Oxide Fuel Cell}
\acrodef{so}[SO]{Solid-Oxide}
\acrodef{soc}[SOC]{Solid-Oxide Cell}
\acrodef{el}[Ely]{Electrolyzer}
\acrodef{fc}[FC]{Fuel Cell}
\acrodef{res}[RES]{Renewable Energy Source}
\acrodefplural{res}[RESs]{Renewable Energy Sources}
\acrodef{ess}[ESS]{Energy Storage System}
\acrodef{pth}[P2H]{Power to Hydrogen}
\acrodef{pv}[PV]{Photovoltaic}
\acrodef{soh}[SoH]{Hydrogen Level in the Tank}
\acrodef{minlp}[MINLP]{Mixed Integer Non Linear Problem}
\acrodef{milp}[MILP]{Mixed Integer Linear Problem}
\acrodef{sos2}[SOS2]{Special Ordered Set of type 2}
\acrodef{mpc}[MPC]{Model Predictive Control}
\acrodef{ec}[EC]{European Commission}
\acrodef{rec}[REC]{Renewable Energy Community}
\acrodef{dmc}[DMC]{Discrete Markov Chain}
\acrodef{igsc}[IGSC]{Instantaneous Growing Stream Clustering}
\acrodef{mise}[MISE]{Italian Ministry of Economic Development}
\acrodef{gams}[GAMS]{General Algebraic Modeling System}

\IEEEoverridecommandlockouts
\newcommand\copyrighttext{%
  \footnotesize
  \centering\copyright~2022 IEEE. Personal use of this material is permitted. Permission from IEEE must be obtained for all other uses, in any current or future media, including reprinting/republishing this material for advertising or promotional purposes, creating new collective works, for resale or redistribution to servers or lists, or reuse of any copyrighted component of this work in other works. \\ DOI: 10.1109/PMAPS53380.2022.9810605}

\newcommand\copyrightnotice{%
\begin{tikzpicture}[remember picture,overlay]
\node[anchor=south,yshift=0pt] at (current page.south) {\setlength{\fboxrule}{0pt}\fbox{\parbox{\dimexpr\textwidth-\fboxsep-\fboxrule\relax}{\copyrighttext}}};
\end{tikzpicture}%
}

\begin{document}

\title{Optimal Management of Renewable Generation and Uncertain Demand with Reverse Fuel Cells by Stochastic Model Predictive Control\thanks{This work was carried out in the framework of the grant PRIN-2017K4JZEE ``Planning and flexible operation of micro-grids with generation, storage and demand control as a support to sustainable and efficient electrical power systems: regulatory aspects, modelling and experimental validation" financed by the Italian Ministry for Education, University and Research.}}

\author{%
  \IEEEauthorblockN{F. Conte}
  \IEEEauthorblockA{Campus Bio-Medico University of Rome\\
    Faculty of Engineering\\
    Via Alvaro del Portillo, 21\\
    I-00128, Roma, Italy\\
    f.conte@unicampus.it}
          \and
     \IEEEauthorblockN{G. Mosaico, G. Natrella, M. Saviozzi}
 \IEEEauthorblockA{University of Genoa\\
    DITEN\\
    Via all'Opera Pia 11 A\\
    I-16145, Genova, Italy\\
    matteo.saviozzi@unige.it}
    \and
\IEEEauthorblockN{F. R. Bianchi}
 \IEEEauthorblockA{University of Genoa\\
    DICCA\\
    Via Montallegro 1 \\
    I-16100, Genova, Italy\\
    fiammettarita.bianchi@edu.unige.it}
    }

\IEEEaftertitletext{\copyrightnotice\vspace{1.1\baselineskip}}

\maketitle
\IEEEpubidadjcol
\begin{abstract}
This paper proposes a control strategy for a Reverse Fuel Cell used to manage a Renewable Energy Community. A two-stage scenario-based Model Predictive Control algorithm is designed to define the best economic strategy to be followed during operation. Renewable energy generation and users' demand are forecasted by a suitably defined Discrete Markov Chain based method.
The control algorithm is able to take into account the uncertainties of forecasts and the nonlinear behaviour of the Reversible Fuel Cell.
The performance of proposed approach is tested on a Renewable Energy Community composed by an aggregation of industrial buildings equipped with PV.
\end{abstract}



\section{Introduction}\label{sec:Introduction}
During last years, the wide spread of \acp{res} has led academic and industrial research to investigate methodologies and technologies which allow a better use of renewable generation to supply energy systems. In literature, different techniques have been studied to manage \ac{res} generation and to optimize the functioning. \acp{res}, such as wind and solar, are variable and hard to predict, therefore many stochastic algorithms have been developed to optimally manage the uncertainties in their forecasts.

The integration of \acp{ess} is necessary to deal with \ac{res} forecasting errors and uncertainties in power demand, and to obtain power system flexibility, namely the ability of the system generators to react to unexpected changes in load or system component performance \cite{IEA2008}. Electrochemical \acp{ess}, such as, batteries, have been widely studied and many works on batteries management can be found in literature \cite{Conte:2019}. A valid and environmental-friendly alternative to batteries are \ac{pth} systems in which possible generation surplus is transformed into hydrogen by an \ac{el} and stored in a tank \cite{gallo2016}. The same hydrogen can be eventually converted back into electrical power by a \ac{fc} when power demand exceeds power generation. Depending on the technology, the \ac{el} and the \ac{fc} can be two different devices or a single reversible unit, called as \ac{rfc}, working in \ac{fc} or \ac{el} mode alternatively \cite{Soloveichik2014}. 

\ac{rfc} dynamic behaviour is more complex than the one of batteries. Its efficiency depends on the working point according to a nonlinear law and the switching between the two modes cannot be executed close to instantaneously as for batteries, but a transition step has to be considered. In this framework, the objective of this study is to develop an optimal control algorithm able to take into account the uncertainties of \ac{res} generation and users demand. In order to satisfy this requirement, we provide a detailed local model for \ac{rfc} operation efficiency and a \ac{dmc}-based forecast algorithm \cite{PMAPS2020} for both load and renewable production. Moreover, we use a two-stage scenario-based programming approach \cite{Parisio:2015,Parisio:2016} to deal with nonlinearities and forecast uncertainties. The result is a stochastic \ac{mpc} algorithm which optimizes the economic operation of a \ac{rec}. The integration of these different approaches and technologies in the \ac{rec} framework is not so common in previous works, which focused on an effective formulation of control system algorithms, not using a \ac{rfc} detailed model, or on \ac{rfc} behaviour without optimizing global operation \cite{Zhang:2019, Mastropasqua:2020}

A \ac{rec} is a legal entity introduced by the \ac{ec} through the Clean energy for all Europeans package. In particular, the \ac{ec} issued two directives IEM \cite{IEM} and RED II \cite{REDII} aiming at improving the uptake of energy communities, at making easier for citizens to integrate efficiently in the electricity system as active participants, and at strengthening the role of \ac{res} self-consumers and \ac{rec}.

In this paper, we consider as case study a \ac{rec} composed of a \ac{pv} power plant, a \ac{rfc} unit based on \ac{soc} technology to cope with uncertainties in the \ac{res} generation and power demand, and an aggregation of industrial warehouses operating as consumers. The manager of the \ac{rec} administrate both \ac{res} generation and \ac{rfc} operation, according to Italian transposition of European directives IEM and RED II.

The paper is organized as follows: section~\ref{sec:SystemArchitecture} introduces the system model, section~\ref{sec:Algorithm} provides the control strategy, the case study is described in section~\ref{sec:StudyCase}, section~\ref{sec:Results} reports simulation results and the conclusions are reported in section~\ref{sec:Conclusions}.

\section{System Model}\label{sec:SystemArchitecture}

\begin{figure}[t]
	\centering
    \includegraphics[width=1\columnwidth]{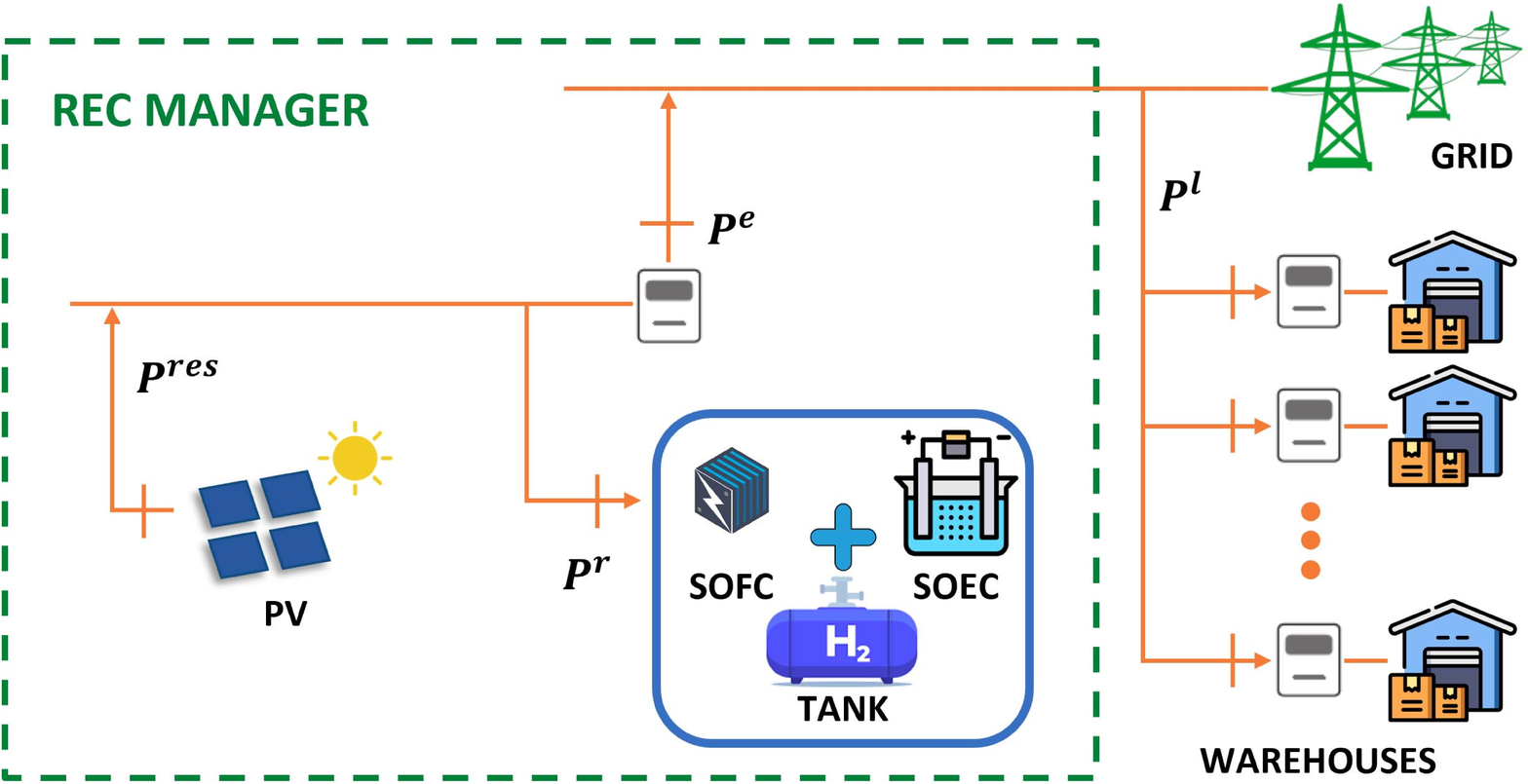}
    	\caption{System Architecture.}
	\label{fig:sys_arc}
\end{figure}

The schematic architecture of the considered system is reported in Figure~\ref{fig:sys_arc}. The \ac{rec} is composed of a \ac{pv} power plant and a \ac{rfc} serving an aggregation of industrial warehouses. According to EU directives, the \ac{rec} can supply the consumers' power demand and also export power to the grid. The industrial aggregation absorbs power from the \ac{pv} plant and eventually from the grid, when the \ac{res} generation does not meet its demand. In the following, the models adopted for each system component are provided. In all, $t$ indicates the discrete-time with a sampling time $\Delta$=\SI{15}{\minute} and reported powers are considered as mean values within the sampling interval.

\subsection{Connection with the Main Grid}
During the quarter hour $t$, the manager of the \ac{rec} can export power $P^e_t$ to the main grid. Therefore, it results that
\begin{equation}
    0 \leq P^e_t \leq P_{max}^e, \label{eq: pe_lim}
\end{equation}
where $P_{max}^e$ is the maximum power exportable from \ac{rec}.

\subsection{\ac{rfc}}
The \ac{rfc} absorbs power $P^{el}_t$ to feed a tank with hydrogen as \ac{soec} or generates power $P^f_t$ by consuming stored hydrogen as \ac{sofc}. \acp{soc} usually work at low pressures and high temperatures reducing crack formation probability and allowing high efficiencies. On the other hand a long start-up is required to reach requested operative temperatures with a time range depending on system size. \ac{soec} and \ac{sofc} have different nonlinear efficiencies. Furthermore when switching from a mode of functioning to the other, the \ac{rfc} does not produce nor consume hydrogen, but it demands power to maintain constant its temperature. We indicate with $\Tilde{P}^{el}$, the power supplied to the \ac{rfc} when switching from \ac{sofc} to \ac{soec}, and with $\Tilde{P}^f$, the power supplied to the \ac{rfc}, when switching from \ac{soec} to \ac{sofc}.
Due to slow thermal response, in this application we have decided to set the \ac{rfc} always on. Equations describing the \ac{rfc} functioning and hydrogen tank managing are reported below:
\begin{align}
    &P^r_t = P_t^{el}\delta_t^{el} - P_t^f\delta_t^f + \Tilde{P}^{el}\Tilde{\delta}_t^{el}  + \Tilde{P}^{f}\Tilde{\delta}_t^{f}, \label{eq: pr} \\
    &P_{min}^f\delta_t^f \leq P^f_t \leq P_{max}^f\delta_t^f, \label{eq: pfc_lim} \\
    &P_{min}^{el}\delta_t^{el} \leq P^{el}_t \leq P_{max}^{el}\delta_t^{el}, \label{eq: pel_lim} \\
    &\delta_t^{el} + \delta_t^f + \Tilde{\delta}_t^{el}  + \Tilde{\delta}_t^{f} = 1,\label{eq: deltas}\\
    &H_{t+1} = H_{t} + \frac{\Delta}{E^h} \left( \phi^{el}_t- \phi^f_t \right), \label{eq:soh} \\
    &H_{min} \leq H_t \leq H_{max}, \label{eq: soh_lim}
\end{align}
where: $P^r_t$ is the power exchanged by the \ac{rfc}, positive when absorbing, negative when generating; $P^{el}_{max}$, $P^{el}_{min}$, $P^f_{max}$ and $P^f_{min}$ are, in the following order, maximum and minimum power of \ac{soec} and \ac{sofc}; $H_t$ [p.u.] is the \ac{soh}; $H_{min}$ and $H_{max}$ are minimum and maximum \acp{soh}; $E^h$ [Wh] is the tank capacity, defined according to the transformation \SI{1}{\mega\watt\hour} = \SI{30}{\kilo\gram}; $\phi^{el}$ and $\phi^f$ are the power exchanged by the \ac{rfc} with the tank, respectively in \ac{soec} and \ac{sofc} mode. $\phi^{el}$ and $\phi^f$ are nonlinear functions of $P^{el}_t$ and $P^f_t$, respectively.

$\delta^{el}_t$, $\delta_t^f$, $\Tilde{\delta}_t^{el}$ and $\Tilde{\delta}_t^{f}$ are binary variables representing the operating mode of the \ac{rfc}, respectively: \ac{soec} mode, \ac{sofc} mode, \ac{tsoec} mode and \ac{tsofc} mode.
In order to switch from \ac{soec} to \ac{sofc}, the \ac{rfc} has to operate in \ac{tsofc} mode before operating in \ac{sofc} mode, similarly in order to switch from \ac{sofc} to \ac{soec}, the \ac{rfc} has to operate in \ac{tsoec} mode first. Furthermore, when \ac{tsofc} or in \ac{tsoec} mode occurs, the \ac{rfc} must operate as a \ac{sofc} or as a \ac{soec}, respectively, in the following time interval. Finally in order to curtail mechanical and thermal stress, the number of switches between the operating modes should be limited. All of these conditions are modeled with the following mixed integer constraints:
\begin{align}
    &\delta_t^{el} + \delta_{t+1}^f \leq 1, \qquad \delta_t^{el} + \Tilde{\delta}_{t+1}^{el} \leq 1, \label{eq: delta_el}\\
    &\delta_t^f + \delta_{t+1}^{el} \leq 1, \qquad \delta_t^f + \Tilde{\delta}_{t+1}^f \leq 1,\label{eq: delta_f}\\
    &\Tilde{\delta}_t^{el} - \delta_{t+1}^{el} \leq 0, \qquad \Tilde{\delta}_t^f - \delta_{t+1}^f \leq 0, \label{eq: delta_tf}\\
    & \sum_{j=0}^{M} \left( \Tilde{\delta}_{t+j}^{el} \right) \leq a, \qquad \sum_{j=0}^{M} \left( \Tilde{\delta}_{t+j}^{f} \right) \leq a. \label{eq: delta_t_lim}
\end{align}
Constraint \eqref{eq: delta_t_lim} is introduced to limit mechanical and thermal stress; $\alpha$ is the maximum number of switches allowed in $M$ time-steps from the current time $t$.


\subsection{\ac{rec}}
According to Italian transposition of REM and RED II \cite{it-trans}, a \ac{rec} is paid for the energy self-consumed between the members of the community and for the energy sold to the grid. The self-consumed energy $\Delta\cdot P^{ac}_t$ is defined as the minimum between the energy exported by the manager and the one consumed by the members of the \ac{rec}:
\begin{equation}
    P^{ac}_t = \min ( P^e_t, P^l_t) \label{eq: pac_lim}
\end{equation}
where $P^l_t$ is the power demand at time $t$ of the warehouses aggregate. 

Since according to the manager economic return \eqref{eq:manager_income}, introduced below, both $P^{ac}_t$ and 
$P^e_t$ will be maximized, \eqref{eq: pac_lim} can be rewritten by the following inequalities: 
\begin{align}
    P^{ac}_t &\geq 0, \\
    P^{ac}_t &\leq P^e_t,  \label{eq: pac_lim1}\\
    P^{ac}_t &\leq P^l_t.\label{eq: pac_lim2}
\end{align}

\subsection{Power Balance, Operational Costs and Available Data}
During every quarter hour $t$, the following power balance has to be matched:
\begin{equation}
    P^{res}_t = P^r_t + P^e_t, \label{eq:powerbalance}
\end{equation}
and the manager economic return is:
\begin{equation}
    J_t = \Delta\left(\left(c^m + c^r\right) P^{ac}_t + c^e_t P^e_t \right) \label{eq:manager_income}
\end{equation}
where $P_t^{res}$ is the power generated by \ac{res}; $c^e_t$ is the energy sell-back price; $c^m$ is an incentive bestowed by the \ac{mise} and $c^r$ is the restitution of grid charges since $P^{ac}_t$ does not burden on the grid \cite{it-trans}. The objective of the paper is to maximize the \ac{rec} manager economic return, assuming that at time step $t$, given a time horizon $T$, the following data are available:
\begin{itemize}
    \item a set of $S$ scenarios each one containing a forecast profile of \ac{res} generation $\{\hat{P}^{res}_{t+k}\left(s\right)\}_{k=0}^{k=T-1}$, a forecast profile of load $\{\hat{P}^{l}_{t+k}\left(s\right)\}_{k=0}^{k=T-1}$, and an associated confidence probability $\pi_t(s)$ associated at the mentioned scenario $s = 1 \dots S$;
    \item the current \ac{soh} $H_t$;
    \item all energy prices from time $t$ to time $t+T-1$.
\end{itemize}

\section{Optimal Management}\label{sec:Algorithm}
In this section we propose the optimal management algorithm, which decides a control action at each time-step $t$, given the data above reported. At quarter hour $t$, we will indicate with $k=0,1,\ldots,T-1$ the time sequence $t,t+1,\ldots,t+T-1$. 

In the following, we first introduce the method adopted to obtain load and \ac{res} forecasts,and then we formulate a \ac{milp}, finally used by an \ac{mpc} controller to perform the optimal management.

\subsection{Load and \ac{res} Generation Forecasts}\label{ssec:forecasts}
\modifPUMA{In \cite{PMAPS2020}, a methodology named \ac{igsc} has been proposed to model time series of interest with a \ac{dmc} through an adaptive online algorithm with minimal computational efforts. The constructed \ac{dmc} can then be used to sample possible future scenarios (forecasts) given the current actual state of the \ac{dmc}.}

\modifPUMA{The states dwell in the same space of the measurements (e.g. in the active power-time plane) and are characterized by the mean of the measurements that happened to be closest to that state. For each state also the number of measurements, their variance, and covariance between the variables are kept in memory.}

\modifPUMA{The algorithm presents just one parameter: $\tau$, a positive real number, which regulates how different from current states must be a new measurement in order to add a new state to the \ac{dmc} (the optimal choice of $\tau$ has been investigated in \cite{AEIT2020}). The algorithm is detailed in the following for the case of two variables.}

\modifPUMA{Let $m(x)=[m_1(x),m_2(x)]$ the mean vector of a state $x$, $\sigma(x)=[\sigma_1(x),\sigma_2(x)]$ the variances of a state $x$, $\rho(x)$ the covariance of a state $x$ and $N(x)$ the measurement assigned to a state $x$. Let $O_k=[O_{1,k}, O_{2,k}]$ be the current observation and $x_{k-1}$ the last timestep state. The steps of the proposed algorithm for each incoming measurement are the following:}
\modifPUMA{\begin{enumerate}[leftmargin=*]
\item Matching step: Find the two states closest to $O_k$, respectively $F_k$ and $S_k$;
    \item State Adaptation: create a new state with the same coordinates of $O_k$, connect it with $F_k$ and assign it to $x_k$ (current state) if: 
    \begin{itemize}
        \item $O_k$ is outside the circle of diameter $\overline{F_k S_k}$; 
        \item The Euclidean distance between $O_k$ and $F_k$ is greater than $\tau$;
    \end{itemize}
        otherwise do not create a new state, instead let $x_k = F_k$
    \item Weight adaptation: The state $x_k$, to which $O_k$ has been assigned, is updated. For each variable $i={1,2}$, the mean is updated as follows: 
    \begin{equation}
        m_i(x_{k}):=\frac{m(x_{k})\cdot N(x_k)+O_k}{N(x_k)+1}
    \end{equation}
    With similar formulas the variance of each variable and the covariance between each pair of variables are updated for current state $x_k$:
      \begin{equation}
        \sigma_i(x_k):=\frac{\sigma_i(x_{k})\cdot N(x_k)+(O_k-m_i(x_k))^2}{N(x_k)+1}
    \end{equation}
      \begin{equation}
      \begin{split}
        \rho(x_k):=\frac{1}{N(x_k)+1}\cdot\rho(x_{k})\cdot N(x_k)\qquad\qquad\\+\frac{1}{N(x_k)+1}(O_{1,k}-m_1(x_k))\cdot(O_{2,k}-m_2(x_k))
        \end{split}
    \end{equation}
     Finally, $N(x_k)$ is increased by one.
    \item Edge Adaptation: A link between the past state $x_{k-1}$ and the current state $x_k$ is created, if it does not exist. Moreover, the weights of all the links starting from the past state are changed to reflect the transition probability. 
\end{enumerate}}

\modifPUMA{The resulting \ac{dmc} is used for the simulation of possible future scenarios, given the current state and using the transition probabilities between the states.}

\modifPUMA{In order to add some variability, the simulated values are not exactly equal to the mean of the states, but they are added to a realization of a bivariate normal random variable having as covariance matrix the variance and covariances computed for each state during the weight adaptation steps.}


To be used by the control algorithm, the \ac{dmc} is first trained on historical data. Then, at time $t$, given the actual values of \ac{res} generation $P^{res}_t$ and power demand $P^l_t$, a set of 300 paths of length $T+1$ is generated. The set is then reduced to 10 scenarios with same length and a probability $\pi_t(s)$ is associated to each of them, applying the scenario reduction method introduced in \cite{Growe:2003}.

\subsection{Piecewise linearization of nonlinearities}\label{ssec:piecewise_linearization}
Equation \eqref{eq:soh}, for any $t=k$, is nonlinear, since $\phi^{el}_k$ and $\phi^{el}_k$ are nonlinear functions of $P_k^{el}$ and $P_k^{f}$, \textit{i.e.} $\phi^{el}_k=g^{el}(P_k^{el})$ and $\phi^{el}_k=g^{f}(P_k^{el})$. To cope with this, \ac{sos2} variables are introduced. \ac{sos2} is an ordered set of non-negative variables, of which at most two of them can be non-zero and if two are non-zero they must be contiguous in the ordered set. Given the set $\Lambda=\{\lambda^l\}_{l=1}^{L}$ of length $L$ it has to be:

\begin{align}
    &\sum^{L}_{l=1}\lambda^l = 1, \quad \lambda^l \geq 0 \\
    \mathbf{if} \quad &\lambda^{l'} > 0:
    \begin{cases}
        \lambda^{\bar l} = 0 \\
        \lambda^{l'+1} \geq 0 \\
    \end{cases}
    \bar l \in [1,L] \land \bar l \neq l,l+1 \\
    & l' \in [1,L-1] \nonumber
\end{align}

\ac{sos2} variables are adopted to approximate functions $g^{el}(\cdot)$ and $g^{f}(\cdot)$, as it follows:

\begin{equation}
    P^{\mu}_k = \sum^{L}_{l=1}\lambda^l_k P^{\mu,l}, \quad
    \phi^{\mu}_k = \sum^{L}_{l=1}\lambda^l_k \phi^{\mu,l} \label{eq:sos_val}
\end{equation}
where $\mu=el,f$, $P^{\mu,l}$ are the independent variable breakpoints, $\phi^{\mu,l}$ are the value of functions at the breakpoints (intercepts: $\phi^{\mu,l}=g^{\mu}(P^{\mu,l})$), and $\lambda_k^l$ are the \ac{sos2} variables.


\subsection{Two-stage Stochastic Optimization Problem}
In a two-stage stochastic programming approach the sum of two cost functions is minimized, the first-stage one refers to the actual objective of the optimization, the second-stage one is suitably defined to minimize the expected violation of constraints that involve random variables. Specifically, such uncertain constraints are relaxed by introducing positive auxiliary variables, called \textit{recourse variables}. The value of these variables is then minimized according to second-stage cost function, that is suitably defined to take into account the probability distributions of the random variables. For details, the reader is referred to \cite{Parisio:2015,Parisio:2016}. 


According to \eqref{eq:manager_income}, first-stage cost function is:
\begin{equation}
    J^{fs} = -\sum_{k=0}^{T-1} \Delta\left(\left(c^m + c^r\right) P^{ac}_k + c^e_k P^e_k \right)  \label{eq:fist-stage-cost-func}
\end{equation}
and constraints are \eqref{eq: pe_lim}--\eqref{eq: delta_t_lim}, \eqref{eq: pac_lim2}--\eqref{eq:powerbalance}, \eqref{eq:sos_val}. Constraints that involve random variables are 
\eqref{eq: pac_lim2} and \eqref{eq:powerbalance}. Therefore, first we rewrite \eqref{eq: pac_lim2} as equality constraint, $\forall k \in [0,T-1]$:
\begin{align}
    P^{ac}_k+\gamma_k=P^l_k, \quad \gamma_k \geq 0, \label{eq:Paclim_eq}
\end{align}
where $\gamma_k$ are slack positive variables; then, we rewrite \eqref{eq:Paclim_eq} and \eqref{eq:powerbalance}, as second-stage constraints: $\forall k \in [0,T-1]$ and $\forall s \in [1,S]$,
\begin{align}
    &\xi_k^+(s) \geq {P}^{ac}_k + \gamma_k -\hat{P}^l_k(s), \label{eq:xiplus1}\\ 
    &\xi_k^+(s) \geq 0, \label{eq:xiplus2}\\
    &\xi_k^-(s) \geq \hat{P}^l_k(s) -{P}^{ac}_k- \gamma_k, \label{eq:ximinus1}\\
    &\xi_k^-(s) \geq 0, \label{eq:ximinus2}\\
    &\gamma_k\geq 0, \label{eq:gamma_positive}\\
    &\chi_k^+(s) \geq {P}^r_k+{P}^e_k-\hat{P}^{res}_k(s), \label{eq:chiplus1}\\ 
    &\chi_k^+(s) \geq 0, \label{eq:chiplus2}\\
    &\chi_k^-(s) \geq \hat{P}^{res}_k(s)-{P}^r_k-{P}^e_k, \label{eq:chiminus1}\\
    &\chi_k^-(s) \geq 0. \label{eq:chiminus2}
\end{align}
where $\chi_k^+(s)$, $\chi_k^-(s)$, $\xi_k^+(s)$ and $\xi_k^-(s)$ are the recourse variables; finally, the second-stage cost function is defined as follows:
\begin{align}
    &J_t^{ss} = \sum_{k=0}^{T-1}  \sum_{s=1}^S\pi_t(s)\left[\left(\chi_k^{+}(s)+ \xi_k^{+}(s)\right) \omega^+ \right. \nonumber \\
    & \qquad \qquad \qquad \qquad \qquad \left. +\left(\chi_k^{-}(s)+ \xi_k^{-}(s)\right)\omega^- \right], \label{eq:costfun2}
\end{align}
where $\omega^+$ and $\omega^-$ are penalty weights. In \eqref{eq:costfun2} we can observe how, using the probability distribution $\pi_t(s)$ of scenarios, recursive variables are minimized as higher is the probability that the corresponding scenarios occur.

To conclude, the two-stage stochastic optimization problem is formulated by the following \ac{milp}:
\begin{equation}
    \begin{aligned}
        &\min_{\{X_k\}_{k=0}^{T-1}} \left(J_t^{fs}+J_t^{ss}\right)\\
        &X_k=\left[{P}^e_k, {P}^{el}_k, {P}^f_k, \delta^{el}_k, \delta^f_k,\Tilde{\delta}^{el}_k, \Tilde{\delta}^f_k \right]^\top \label{eq:cost_function}
    \end{aligned}
\end{equation}
subject to:
\eqref{eq: pe_lim}--\eqref{eq: delta_t_lim}, \eqref{eq: pac_lim1}, \eqref{eq:sos_val}, \eqref{eq:xiplus1}--\eqref{eq:chiminus2}.

\subsection{Control Algorithm}
The algorithm consists in solving \eqref{eq:cost_function} at each time-step $t$. Then, the \textit{receding horizon} is adopted: to apply just the values calculated for the instant $k=0$ to the controlled variable, move to the successive time-step $t+1$, and repeat the same procedure. Forecasting errors are compensated outside the optimization problem with the aim of maximizing the income of the \ac{rec} manager. Therefore, if the actual \ac{res} generation is higher than the expected one, $P^e_t$ is increased to exploit the surplus in generation and increase the income; instead, if it is lower, $P^r_t$ is increased, when possible according to \eqref{eq: pe_lim}--\eqref{eq: soh_lim}, otherwise $P^e_t$ is decreased. Changes in $P^r_t$ and $P^e_t$ are made to keep always satisfied \eqref{eq:powerbalance}.



\section{Case Study}\label{sec:StudyCase}
The considered case study is an aggregation of industrial warehouses equipped with \ac{pv} generators. The industrial aggregation is supposed to undergo Italian regulation \cite{it-trans}. Data of power demand and \ac{pv} generation are taken form \cite{data}. In particular, power demand was taken as it was, while \ac{pv} generation was increased in nominal power with respect to the values reported in \cite{data}. Values of nominal power demand and nominal \ac{pv} generation are listed in Table~\ref{tab:parms}.

\begin{figure}[t]
	\centering
    \includegraphics[width=1\columnwidth]{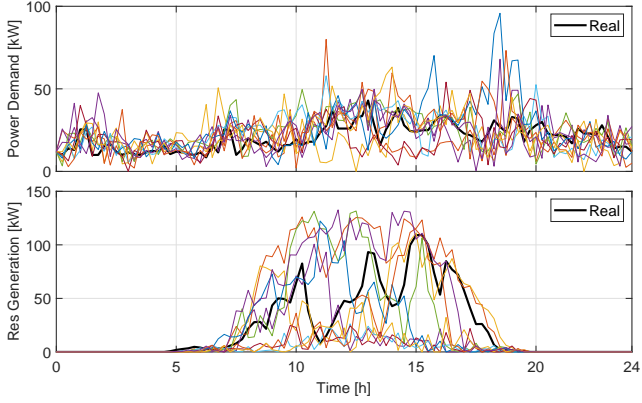}
    	\caption{Example of 10 scenarios generated at midnight of the third day of simulations for warehouses aggregated demand (top) and \ac{res} generation (bottom).}
	\label{fig:scenari}
\end{figure}

To get forecasts, the method described in Section~\ref{ssec:forecasts} is applied with $\tau=0.1$. The \ac{dmc} is trained with historical data in \cite{data}. Figure~\ref{fig:scenari} shows an example of 10 scenarios generated for \ac{res} production and power demand at the midnight of the third day of simulation in comparison with their real values.

The developed model correlates the power, in terms of hydrogen production (\ac{soec}) or consumption (\ac{sofc}), with the number $N_c$ of \acp{soc} constituting the \ac{rfc}, the operating temperature $\theta$ of the \acp{soc}, and the electrical powers used $P^{el}$ and produced $P^f$.  The investigated relations are:

\begin{align}
     &\eta^{el} = \frac{N^{el}_{H_2}LHV_{H_2}}{P^{el}} = \frac{\phi^{el}}{P^{el}} \label{eq:etael}, \\
     &\eta^f = \frac{P^f}{N^f_{H_2}LHV_{H_2}}=\frac{P^f}{\phi^f}, \label{eq:etaf}
\end{align} 
where $\eta^{el}$ and $\eta^f$ are efficiencies in \ac{soec} and \ac{sofc} mode, $N^{el}_{H_2}$ and $N^f_{H_2}$ the flow of hydrogen produced in \ac{soec} mode and consumed in \ac{sofc} mode, $LHV_{H_2}$ the Lower Heating Value of hydrogen, $P^{el}$ the power that has to be given in \ac{soec} mode to produce $N^{el}_{H_2}$ ensuring isothermic operation, warming up the reagents and compressing the produced hydrogen, finally $P^f$ the electric power produced by the \ac{sofc} deducted by the power to warm up the reagents. \ac{soc} working conditions are fixed in order to optimise the operation basing on authors' previous work \cite{Bianchi:2021}.

\begin{figure}[t]
	\centering
    \includegraphics[width=1\columnwidth]{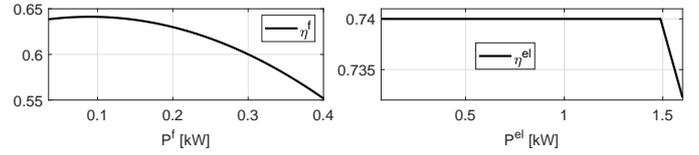}
    	\caption{Single cell efficiencies in \ac{sofc} (left) and \ac{soec} (right) modes, with $\theta=\SI{1123}{\K}$. }
	\label{fig:eta}
\end{figure}

In the following equations we provide the numerical expressions of $\eta^{el}$ and $\eta^f$. 
\begin{align}
     \eta^{el} =& \label{eq:etael2}
    \begin{cases}
        0.74 \quad &\frac{P^{el}}{N_c} \leq \left(\alpha_1\theta^2 - \alpha_2\theta + \alpha_3\right)\\
        \beta_1\frac{\theta P^{el}}{N_{c}} - \beta_2\frac{P^{el}}{N_{c}} + \beta_3\theta \quad &\frac{P^{el}}{N_c} > \left(\alpha_1\theta^2 - \alpha_2\theta + \alpha_3\right)\\
    \end{cases}\\
    \eta^f = &\num{8.06e-3}\frac{\theta P^f}{N_{c}} - 8.89\frac{P^f}{N_{c}} \nonumber + \num{1.85e-2}\theta\\
    &- \num{9.29e-1}\left(\frac{P^f}{N_{c}}\right)^2 - \num{8.95e-6}\theta^2  - 8.88 \label{eq:etaf2}
\end{align}
where $\alpha_1 = \num{4.87e-4}$, $\alpha_2 = \num{9.46e-2}$, $\alpha_3 = 46.34$, $\beta_1 = \num{2.32e-4}$, $\beta_2 = 0.33$ and $\beta_3 = \num{7.7e-4}$. The profiles of the two efficiencies, for the constant operating temperature $\theta=\SI{1123}{\K}$ are are shown in Figure~\ref{fig:eta}.

From \eqref{eq:etael}--\eqref{eq:etael2} and \eqref{eq:etaf2}, it is possible to express functions $\phi^{el}=g^{el}(P^{el})$ and $\phi^f=g^f(P^f)$, then used by the control algorithm as described in Section~\ref{ssec:piecewise_linearization}.

The number of \acp{soc} composing the \ac{rfc} was estimated by setting the upper power in \ac{soec} mode equal to the maximum value of power demand. The maximum number of switches between operating modes is set to 3 in \SI{4}{\hour}.

\begin{table}[t]
		\centering
		\caption{Study Case Parameters}
		\label{tab:parms}
		{\small
		\begin{tabular}{lcc}
		\toprule
		Parameter & Symbol  & Value\\
		\midrule
        Main Grid Connection Nominal Power & $P^g_{max}$ & \SI{340}{\kilo\watt} \\
        Storage Capacity & $E^h$ & \SI{400}{\kilo\watt\hour} \\
        \ac{rfc} Number of \acp{soc} & $N_{c}$ & 100\\ 
        Minimum Power in \ac{soec} mode & $P^{el}_{min}$ & \SI{7.2}{\kilo\watt}\\
        Maximum Power in \ac{soec} mode & $P^{el}_{max}$ & \SI{160}{\kilo\watt}\\
        Minimum Power in \ac{sofc} mode & $P^f_{min}$ & \SI{3.5}{\kilo\watt}\\
        Maximum Power in \ac{sofc} mode & $P^f_{max}$ & \SI{40}{\kilo\watt}\\
        \ac{rfc} Power Demand in \ac{tsoec} mode & $\Tilde{P}^{el}$ & \SI{2.6}{\kilo\watt}\\
        \ac{rfc} Power Demand in \ac{tsofc} mode & $\Tilde{P}^f$ & \SI{1.3}{\kilo\watt}\\
        Nominal Power Demand & $-$ & \SI{163}{\kilo\watt} \\
        \ac{pv} Plant Nominal Power & $-$ & \SI{150}{\kilo\watt} \\
        \ac{mise} Incentive & $c^m$ & \SI{0.11}{\eur/\kilo\watt\hour}\\
        Restitution of Grid Charges & $c^r$ & \SI{0.009}{\eur/\kilo\watt\hour}\\
		\bottomrule
		\end{tabular}
        }
\end{table}

\section{Simulation Results}\label{sec:Results}
\begin{figure}[t]
	\centering
    \includegraphics[width=1\columnwidth]{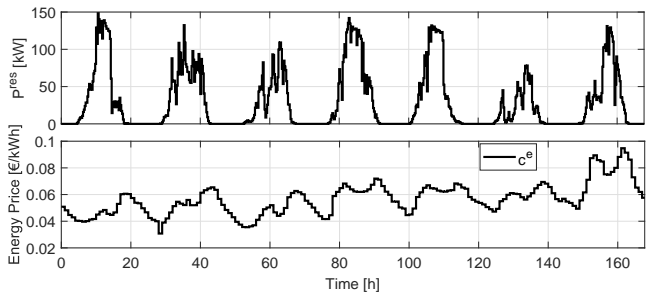}
    	\caption{\ac{res} generation (top) and energy price (bottom) profiles.}
	\label{fig:res-pun}
\end{figure}

To test the performance of the proposed control algorithm, one week has been simulated, using the data from the midnight of June 4th to the midnight of June 12th, 2017. The values adopted for $c^c$, $c^m$ and $c^r$ are listed in Table~\ref{tab:parms}, $c^e_t$ is shown in Figure~\ref{fig:res-pun} and it represents to the actual energy clearing market price in Italy. Notice that the energy prices considered in this application are the characteristic values at the beginning of 2021, before the sudden rise of natural gas price in Europe.

\begin{figure}[t]
	\centering
    \includegraphics[width=1\columnwidth]{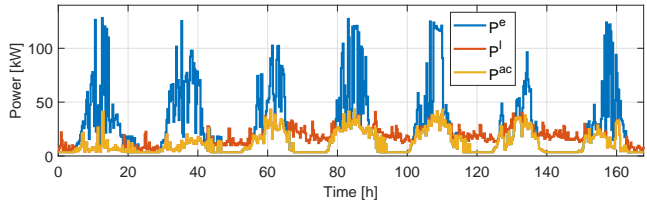}
    	\caption{\ac{rec} power profiles.}
	\label{fig:rec}
\end{figure}

\begin{figure}[t]
	\centering
    \includegraphics[width=1\columnwidth]{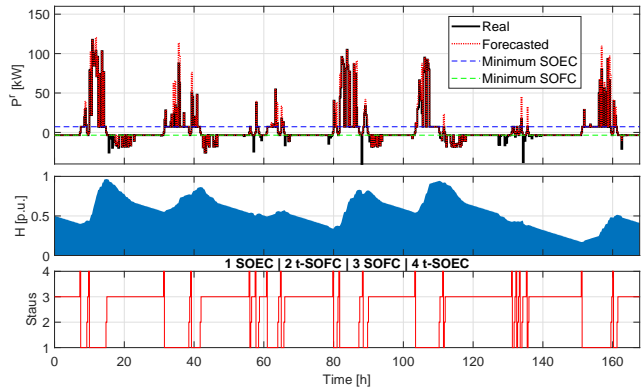}
    	\caption{\ac{rfc} exchanged power (top), \ac{soh} (mid) and operating mode (bottom).}
	\label{fig:rfc}
\end{figure}

The following control parameters have been set as: $T=\SI{24}{\hour}$, $S=10$ scenarios, $H_{max}=\SI{1}{p.u.}$ and $H_{min}=\SI{0}{p.u.}$. The values of penalty weights $\omega^+$ and $\omega^-$ were set to 1, about ten times higher than the values of the other prices.

The control algorithm has been implemented in MATLAB, integrated with \ac{gams} to write the optimization problem, solved by CPLEX solver.

Figures~\ref{fig:rec}--\ref{fig:rfc} show an example of  obtainable simulation results. In particular in Figure~\ref{fig:rfc} we can observe how the \ac{rfc} is managed: during the nights, the algorithm decides to discharge the hydrogen tank and to never switch between two modes since any transition represents an additional load that can not be satisfied in absence of \ac{res} generation. Furthermore we can observe that the majority of switches is obtained during the days with less \ac{res} generation.

We finally remark that the total earning obtained within the week operations has resulted to be equal to \SI{2044.45}{\eur}, against \SI{1955.58}{\eur} obtained in the same conditions without the \ac{rfc}.

\section{Conclusions}\label{sec:Conclusions}
In this paper a control strategy for a \ac{rfc} used to manage a \ac{rec} is proposed. A two-stage scenario-based \ac{mpc} algorithm has been designed to decide the best economic strategy to be followed during operations. Such an algorithm uses a suitably defined \ac{dmc} based method to forecast consumers demand and renewable generation, and a nonlinear model the of \ac{rfc} efficiency derived from a physical based model at local level. The algorithm has been successfully tested on a \ac{rec} composed by an aggregation of industrial buildings and a \ac{pv} plant.


\end{document}